\documentclass[11pt,a4paper]{article}
\usepackage{amsmath}
\usepackage{amsfonts}
\usepackage{amssymb}
\usepackage{xspace}
          \usepackage{graphics}
\usepackage{graphics}

\newtheorem{theorem}{Theorem}

\renewcommand\th{\theta}
\def\l{\left}  \def\r{\right}
\def\a{\alpha} \def\b{\beta}
\renewcommand\t{\tau}

\newcommand\proof{\noindent {\sc Proof:}\qquad}

\newcommand\LB[1]{\label{#1}} 
\newcommand\BE[2]{\begin{#1} #2 \end{#1}}
\newcommand\ARR[2]{\BE{array}{{#1} #2}}
\newcommand\EQ[2]{\BE{equation}{\LB{#1} #2}}

 \newcommand\EQn[1]{\BE{equation*}{ #1}}

\newcommand\EQAn[2]{\EQn{\ARR{#1}{#2}}}

\newcommand{\ay}{asymptotic\xspace}
\newcommand{\eq}{equation\xspace}
\newcommand\wrt{with respect to \xspace}

\newcommand\sm{\sum_1^\iy}
\newcommand{\st}{\, \big| \,  }
\newcommand\ST{such that }
\newcommand{\lint}{\int\limits}

\newcommand\f{\varphi}

\newcommand\gm{\gamma}

\newcommand \Lm{\Lambda}
\newcommand\Th{\Theta}

\newcommand{\fH}{\mathfrak{H}}

\newcommand{\cH}{\mathcal{H}}
\newcommand{\cL}{\mathcal{L}}

\newcommand\eps{\varepsilon}

\newcommand\dl{\delta}

\newcommand{\bbR}{\blackboard{R}}
\newcommand{\blackboard}[1]{\mathbb#1}
\newcommand{\iy}{\infty}

\def\proof{\noindent {\sc Proof:}\qquad}

\title{Regularity of  the  Gurtin-Pipkin  \eq
}
\author{
{S. A. Ivanov\thanks{St. Petersburg Institute of Terrestrial Magnetism,
Ionosphere and Radio Wave Propagation.
{{\tt sergei.a.ivanov@.mail.ru }}}}
\thanks{
The work was supported by Russian Foundation for Basic Research, RFBR
Project 11-01-00790a and RFBR
Project 11-01-00667a. }
}

\begin{document}
\date{}

\maketitle

\begin{abstract}
We study regularity of the solution $\th$ to
the Gurtin-Pipkin integral-differential \eq of the
first order in time.   In particular, we prove
that the 'perturbation' part, namely, the difference of $\th$
and the solution to the corresponding wave \eq is smoother than $\theta$.

\end{abstract}

\vskip1cm
\section{\LB{intro_etc} Introduction}

In several fields of physics such as heat transfer with finite
propagation speed \cite{GuPip}, systems with thermal memory
\cite{EMV}, viscoelasticity problems \cite{CMD}, and acoustic
waves in composite media \cite{Sham}, the  integro-differential
equations arise. We consider the equation of the first order
in time
\begin{equation}\LB1
\theta_{t}(x,t)=\int_0^t k(t-s) \th_{xx}(x,s)\,ds+ f(x,t) ,\
\ x\in(0,\pi), \ t>0,
\end{equation}
with the Dirichlet boundary conditions and with the initial data
$\th(0,x)=\xi(x)$.

In the case $k(t)=Const=\a^2$  the \eq \eqref1 is, in a fact, an integrated wave \eq. Indeed,
differentiate \eqref1  gives
\EQ{wave_eq}{
\th_{tt}=\a^2\th_{xx}+f_t(x,t), \ \th(x,0)=\xi,\ \th_t(x,0)=0.
}
Thus, the wave \eq is a special case of \eqref1 and we will
compare general regularity results
with the regularity of the solutions to the wave \eq.

\section{Fourier method and the Laplace transform}

First, apply the Fourier method:
we set $\f_n=\sqrt{\frac2\pi}\sin nx$ and expand the solution, the RHS, and the initial
data  in series in $\f_n$
$$
\th(x,t)=\sum_1^\iy \th_n(t)\f_n(x),\
 \xi(x)=\sum_1^\iy \xi_n \f_n(x),\
 \ f(x,t)=\sm f_n(t)\f_n(x).
$$

The components   $\th_n$  satisfy ordinary
integral-differential equations

\begin{equation}\label{thn1}
\dot \theta_n(t)=-n^2\int_0^t k(t-s) \theta_n (s)d s\ +f_n(t). \  t>0,
\qquad \theta_n(0)=\xi_n.
\end{equation}

Note that the solutions to this integral-differential is unique and continuous,
what we can see by $t$-integration of the \eq from 0 to $t$. Indeed, we
obtain a Volterra integral equation:  with $\f=\lint_0^t f(t)\,dt$, we have
$$
\th_n(t)-\xi_n=-n^2\lint_0^td\t\lint _0^\t k(\t-s)\th_n(s)\,ds +\f(t).
$$
Change order of integrations
\EQ{2.1}{
\th_n(t)=-n^2\lint_0^tds \,\th_n(s)\lint _s^t \,d\t \,k(\t-s) +\f(t)+\xi_n
=-n^2\lint_0^tq(t-s)
\,\th_n(s)ds +\f(t)+\xi_n
}
with
$$
q(s)=\lint_0^s k(y)\,dy.
$$
Introduce the scale $\cH_s$, $s$ is real, of the Hilbert spaces
$\cH_s=Dom(A^{s/2})$,
where the operator $A$ is  $-d^2/dx^2$ with the Dirichlet
boundary conditions at $0$ and at $\pi$.
A space $\cH_s$ is a subspace of the Sobolev space $H^s$ and
may be described in terms of the Fourier coefficients. Let
the space $l_s$ be the space of sequences $\{c_n\}$ \ST
$$
 \sm|c_n|^2n^{2s}<\iy.
$$
Then
$$
\cH_s=\Big\{u(x)=\sm u_n\f_n(x)\st \{u_n\}\in l_s\Big\}.
$$
Consider also the space $\fH_{s,\eps}$ of functions   $g(x,t)=\sm f_n(t)\f_n(x)$
with the norm
$$
\|g\|^2_{\fH_{s,\eps}}=\sm n^{2s}\|e^{-2\eps t}f_n\|_{L^2(0,\iy)}^2.
$$
\BE{definition}{
The function $\th(x,t)=\sm \th_n(t)\f_n(x)$ is a solution to \eqref1 in
$\fH_{s,\eps}$ if the functions
$\th_n$, satisfy the integral \eq \eqref{2.1} and $\th\in \fH_{s,\eps}$ with $s\in\bbR$.}

Let $H^2_\eps$ denote the
Hardy space in the right half  plane $\Re z>\eps$. The Paley-Wiener theorem says that
$$
\|F\|^2_{H^2_\eps}=\lint_0^\iy e^{-2\eps t}|f(t)|^2 dt.
$$

Here and in what follows will denote the Laplace image by the capital characters.
Applying the Laplace Transform to \eqref{thn1}
we find
\EQn{
z\Theta_n(z)-\xi_n=-n^2 K(z)\Theta_n(z)+F_n(z)
}
or
\EQ{GP1Thn}{
\Th_n(z)=\frac{\xi_n+F_n(z)}{z+n^2K(z)}.
}
Denote the denominators in \eqref{GP1Thn}  by $G_n(z)$.
The set $\Lm$ of all zeros of $G_n(z)$
is called the \emph{spectrum} of the equation \eqref1.

Regularity of the Gurtin-Pipkin type \eq
is studied in \cite{P05} for several spatial variables, where under assumption that
$k(t)$ is twice continuously differentiable it was shown,
in particular, that $\Th(x,t)\in C([O,T];L^2(0,T))$.
Regularity of strong solutions has been studied in several works of
V. Vlasov with the coauthors, see, e.g., \cite{RSV} and the Sec.\ref{Strong} below.
The spectrum of the equation
is studied in \cite{IE},\cite{RSV}.


Let us describe the regularity of the solutions to the
 the wave equation \eqref{wave_eq}. Let $Q_T=(0,\pi)\times (0,T)$.
\BE{proposition}{ The solution to the \eqref{wave_eq} satisfy the following estimates:
(i) Let $f=0$. Then the Dalambert solution gives
\EQ{regwe1}{
\|\th\|_{\fH_{s,\eps}}\prec\|\xi\|_{\cH_s},\
\|\partial_t \th\|_{\fH_{s,\eps}}\prec\|\xi\|_{\cH_{s+1}},
}
(ii)
let $\xi=0$. Then, see \cite{Mikh} the  (generalized) solutions satisfy
\EQ{regwe2}{
\|\partial_t \th\|_{L^2(Q_T)}+
\|\partial_x \th\|_{L^2(Q_T)}
\prec\|\xi\|_{\cH_{1}}+\|f_t\|_{L^2(Q_T)}.
}
}

If $k(t)=\a^2e^{-bt}$ (and $K(z)=\a^2/(z+b)$),
then differentiation gives a damped wave \eq
\EQ{wave1}{
\th_{tt}=\a^2\th_{xx}-b\th_t.
}
By $\th^0$ we denote the solution to this \eq with the initial data
\EQ{waveBC}{
\th(\cdot,0)=\xi, \ \th_t(\cdot,0)=0.
}
This will be an unperturbed \eq, see the Sect. \ref{perturb}.

In application,  see, e.g.,\cite{Sham}, the kernels $k(t)$ is a series of exponentials
$$
k(t)=\sum_1^\iy a_ke^{-b_kt}, \ a_k\ge0,\ 0\le b_1<b_2<\dots<b_k<\dots.
$$
We can consider the following smoothness conditions

\EQ{C0}{
\sum_1^\iy \frac {a_k}{b_k}<\iy,
}

\EQ{C1}{
\a^2=\sum_1^\iy a_k<\iy,
}
or
\EQ{C2}{
\b=\sum_1^\iy a_k b_k<\iy,
}
or
\EQ{C3}{
\gm=\sum_1^\iy a_k b_k^2<\iy.
}
\BE{remark}{These conditions maybe written as
\EQAn{l}{
k\in  L^1(0,\iy),\\
k\in C[0,\iy),\ k,\, k'\in  L^1(0,\iy),\\
k'\in C[0,\iy),\ k,\, k',\, k''\in  L^1(0,\iy),\\
k''\in C[0,\iy),\ k,\, k',\, k'',\, k'''\in  L^1(0,\iy).
}
}

Write the \ay of $K(z)$.
The Laplace image of $k(t)$ is
$$
K(z)=\sum_1^\iy \frac{a_k}{z+b_k}, \ k(0)=\a=\sm a_k.
$$
Without loss of generality we can
set $\a=1$ if $\a$ is finite.
\BE{proposition}{
Let for a $\dl>0$
\EQn{
|\arg z|<\pi-\dl,
}
Then for large $z$

(i) under \eqref{C0}
\EQn{
K(z)=o(1),
}

(ii) under \eqref{C1}
\EQn{
K(z)=\frac1{z}+o(\frac1{z}),
}

(iii)  under \eqref{C2}
\EQn{
K(z)=\frac1{z}-\frac\b{z^2}+o(\frac1{z^2}),
}

(iv)
 under \eqref{C3}
\EQ{ay2}{
K(z)= \frac{1} z-\frac\b{z^2}+\frac \gm{z^3}+o(\frac1{z^3}).
}
}

The statement of this proposition follows from known results about
Cauchy transform of a measure.

\section{Regularity of the  solution in the spatial variable}
Here we prove the results about the regularity \wrt the x-variable, i.e., in terms of
$\cH_s$ spaces.

\begin{theorem}
{\LB {mainthGP1reg}}
Let  \eqref{C1} be true,  $\{\xi_n\}\in\ell_s$ and $f\in L^2(0,\iy;\cH_s)$.
Then for any $\eps>0$
the solution $\th $ to \eqref1 satisfy
\EQ{starp5}{
 \|\th\|^2_{\fH_{s,\eps}} \prec \|\xi\|_s^2
+\|f \|^2_{\fH_{s,\eps}.}
}
\end{theorem}

\proof

\BE{lemma}{\LB{l1}
The following estimates are fulfilled
\EQ{2.1a}{
|z/G_n(z)|\prec1, \ \Re z>\eps
}
\EQ{2.2}{
\|1/G_n\|_{L^2(\eps-i\iy, \eps+i\iy)}\prec1.
}
}
The lemma implies
by  \eqref{GP1Thn}
$$
\lint_0^\iy |e^{-\eps t}\th_n(t) |^2 \prec |\xi_n|^2+ \|e^{-\eps t}f _n\|^2_{L^2(0,\iy)}
$$
and then  \eqref{starp5}.

Proof of the lemma.
Set for the simplicity $\eps=1$    
Then for $z=1+iy$
and $\gm_k=1+b_k$ we obtain
$$
G_n(z)=(1+iy)+n^2\sm\frac{a_k\gm_k}{\gm_k^2+y^2} -iyn^2
\sm\frac{a_k}{\gm_k^2+y^2}
$$
Therefore
$$
|G_n(z)|^2\ge \l( 1+n^2\frac{a_1\gm_1}{\gm_1^2+y^2} \r)^2
+y^2\l( 1- n^2\sm\frac{a_k}{\gm_k^2+y^2}\r)^2
$$
Setting
$$
s(y)=\sm\frac{a_k}{\gm_k^2+y^2}.
$$
we have
$$
|G_n(z)|^2\succ \l( 1+n^2\frac{1}{1+y^2} \r)^2+y^2\l( 1- n^2s(y)\r)^2.
$$
This gives \eqref{2.1a}.

Divide $[0,\iy)$ into three intervals
$$
I_1=[0,n/2]\ ,I_2=[n/2,3n/2]\ ,I_3=[3n/2,/iy].
$$
Write
\EQn{
\|1/G_n\|^2_{L^2(\eps-i\iy, \eps+i\iy)}
=\lint_{\iy }^\iy \frac{dy}{|G_n(y)|^2}=
2\l[\lint_{0 }^{n/2}+ \lint_{n/2}^{2n}+\lint_{2n}^\iy \r]\frac{dy}{|G_n(y)|^2}
=2(J_1+J_2+J_3).
}

1. Estimates on $I_1=[0,n/2]$.

Evidently, $s(y)$ decreases and then on $[0,n/2]$ we have $s(n/2)<s(y)<s(0)$.
Further, the series
$$
n^2s(y)=\sm \frac{a_k n^2}{\gm_k^2+n^2/4}
$$
has the majorant 4$\sm a_k=4\a=4 $  and the terms of this series approaches to $4a_k$.
Then
$$
n^2s(n/2)\ \to \  4.
$$
Take $n>n_0$ such that $n^2s(n/2)\ge 2$ for $n/2\ge n_0$. We obtain
$n^2s(y)>n^2s(n/2)\ge 2$ and
$$
\l(n^2s(n/2)-1\r)^2 \ge 1.
$$
This gives
$$
|G_n(z)|^2\ge 1+y^2\l(n^2s(n/2)-1\r)^2\ge 1+y^2.
$$

Estimate $J_1$.
\EQ{J1}{
J_1\le
\lint_0^{n/2} \frac1{1+y^2}\prec1.
}

2. Estimate $J_2$.

For $n/2\le y\le 2n$ we have
$$
|G_n(z)|\succ 1+n^2(1-n^2s(y),
$$
Consider the increasing variable $\xi=1-n^2s(y)$. Then
$$
\xi'=-n^2s'=n^2\sm \frac{2a_ky}{(\gm_k^2+y^2)^2}\asymp n^2 \sm \frac{a_kn}{(\gm_k^2+n^2)^2}
$$
$$
\asymp
\frac1n \sm \frac{a_k}{(\gm_k^2/n^2+1)^2}\asymp\frac1n.
$$
Indeed,
$$
\frac{a_1}{(\gm_1^2/n^2+1)^2}\le
\sm \frac{a_k}{(\gm_k^2/n^2+1)^2}\le \sm a_k.
$$

Now for $J_2$ we have "$nd\xi\asymp dy $" and
$$
J_2\succ
\lint_{n/2}^{2n}\frac{dy}{1+ n^2(1-n^2s(y))^2}=
\lint_{\xi(n/2)}^{\xi(2n)}\frac{nd\xi}{1+n^2\xi^2}<\iy.
$$

3. Estimate $J_3$.

For $y\ge 2n$ we have
$$
n^2s(y)\le n^2s(2n)=\sm\frac{a_kn^2}{\gm_k^2+4n^2}\le \frac14 \sm a_k=\frac14.
$$
Now
$$
|G_m(z)|\ge 1+y^2(1-n^2s(y))^2\ge 1+y^2\frac9{16},
$$
and
$$
J_3\prec
\lint_{2n}^{\iy} \frac{dy}{1+y^2}\prec1.
$$
The theorem is proved.

\BE{remark}{ For the case $k(t)=1$, and $f=0$, i.e.,  for the wave \eq we have
$$
\Th_n(z)=\frac{\xi_n}{z+n^2/z}, \ \th_n(t)=\xi_n\cos nt.
$$
We see that $\th_n\notin L^2(0,\iy)$ and
$e^{-\eps t}\th_n\in L^2(0,\iy)$. In this sense Theorem
\ref{mainthGP1reg} is sharp.
}
\BE{theorem}{\LB{cont}
Let  \eqref{C1} is true,  $\{\xi_n\}\in\ell_s$ and $f\in L^2(0,\iy;\cH_s)$.
Then $\th(x,t)$ is an $\cH_s$ valued continuous function:
$$
\|\th(t)-\th(t+t_0) \|_{\fH_{s,\eps})} \ \to \ 0,
$$
as $t_0\to t_0$.

}
Proof.  The solutions $\th_n$ are continuous and the series in $\th_n$ has a majorant.

\section{\LB{perturb}GP as a perturbation to the wave \eq}

Let us find regularity of the 'perturbation' $\th-\th^0$ of the solution to the wave \eq.
Recall that  $\th^0$   is the solution to the problem \eqref{wave1},\eqref{waveBC}.
Let $f(x,t)=0$ for simplicity and set
$$
K_0(z)=\frac1{z+\b}, \ G_n^0(z)=z+n^2K_0(z),\
D_n(z)=\frac1{G_n(z)}- \frac1{G_n^0(z)}.
$$
The solution to \eqref1 has the form
$$
\Th_n(z)=\frac1{G_n^0(z)}\xi_n+D_n(z)\xi_n= \Th_n^0(z)+D_n(z)\xi_n.
$$

\begin{theorem}
{\LB {mainthGP1reg_wave}}
Let \eqref{C3} is true and $f=0$.
Then for $s<9/2$
$$
\|z^sD_n\|_{L^2(i\bbR)}\prec n^{s-1}|\xi_n|.
$$
\end{theorem}

\proof

If $\b=0$ the theorem is trivial: $D_n=0$. Thus, we can assume $\b\ne0$ and
then integrate the functions on the imaginary axis.

\BE{lemma}{\LB{lm1}
For $z=iy$, $y\to \iy$ and $\b\ne 0$ we have
$$
|G_n(z)|^2\asymp|G^0_n(z)|^2\asymp
\frac1{y^4}\l[ \l(y^2(y^2-n^2)\r)^2+n^4  \r]=:Q(y).
$$
}
Proof of the lemma.
\eqref{ay2} implies
$$
G_n(iy)=iy+n^2\l(
\frac1{iy}-\frac\b{y^2}+o\l(\frac1{y^2}\r)
\r)
$$
$$
=\frac i{y^2}\l[
y^3-n^2y+o(1)n^2
\r]
-\frac{n^2}{y^2}[\b+o(1)].
$$
And the same is true for $G_n^0$.
From here
$$
|G_n(y)|^2=
\frac 1{y^4}\l[
y^3-n^2y+o(1)n^2
\r]^2
+\frac{n^4}{y^2}[\b+o(1)]^2\asymp
\frac 1{y^4}\l[
y(y^2-n^2)+o(1)n^2
\r]^2
+\frac{n^4}{y^2}.
$$
Now
$$
\frac{|G_n(y)|^2}{Q(y)}=\frac{\l[
y(y^2-n^2)+o(1)n^2\r]^2+n^4}{ \l(y^2(y^2-n^2)\r)^2+n^4  }
$$
Use the elementary inequalities
$$
(a+ qb)^2+b^2\asymp a^2+b^2
$$
with, say, $q<1/2$.
Setting
$$
a=y(y^2-n^2), \ b=n^2
$$
we complete the proof of the lemma.

Return to the proof of the theorem.
$$
\|z^sD_n\|_{L^2(i\bbR)}^2=\lint_{-\iy}^\iy n^4
\frac { y^{2s} |K(iy)-K_0(iy)|^2}{|iy+n^2K(iy)|^2 |iy+n^2K_0(iy)|^2}dy.
$$
Use Lemma \ref{lm1} to estimate the denominator and
the condition \eqref{C3} to estimate the numerator.
By \eqref{ay2} we have $K(iy)-K_0(iy)=O(1/z^3)$ and obtain
\EQ{gm1}{
\|z^sD_n\|_{L^2(i\bbR)}^2\prec n^4\lint_0^\iy
\frac { y^{2s+2}}{\l[ y^2(y^2-n^2)^2+n^4\r]^2}dy.
}
For $s<9/2$ this integral converges.
The main contribution gives the interval of the length $O(n)$ centered at $y=n$. Estimate it.
Set
$$
J=\lint_{n/2}^{3n/2}
\frac { y^{2s+2}}{\l[ y^2(y^2-n^2)^2+n^4\r]^2}dy.
$$
(this integral enter \eqref{gm1} with the factor $n^4$).
$$
J\asymp n^{2s-2}\lint_{n/2}^{3n/2}\frac{dy}{\l[ (y^2-n^2)^2+n^2\r]^2}
\asymp n^{2s-2}\lint_{n/2}^{3n/2}\frac{dy}{\l[ n^2(y-n)^2+n^2\r]^2}
$$
$$
\asymp n^{2s-6}\lint_{n/2}^{3n/2}\frac{dy}{\l[ (y-n)^2+1\r]^2}
\le n^{2s-6}\lint_{n/2}^{3n/2}\frac{dy}{\l[ t^2+1\r]^2}=Cn^{2s-6}.
$$

Thus, the
interval of  $(n/2,3n/2)$ gives the contribution $n^{2s-2}$.
The theorem is proved.

This theorem give, of course, information about regularity of $\th-\th^0$, for example
\BE{corollary}{\LB{cc1}}
$$
e^{-\eps t}[\th(x,t)-\th_0(x,t)]\in L^2(0,\iy;\cH_{s+1})\cap W_2^1(0,\iy;\cH_{s})
$$


\section{\LB{Strong}Estimate of the derivative \wrt time}

Before now we have considered a weak solution. Here we present results
about a strong solutions, as, e.g., in \cite{RSV}.
By the definition for the strong solution
the equation \eqref1 can be considered as an equality of elements of $L^2$ spaces.

\BE{proposition}{\LB{VV}}{ Let \eqref{C2} is fulfilled,
$f_t\in \fH_{1,\gm}$, $f(x,0)=0$,
and $\xi\in \cH_2$. Then for any $\eps>\gm$
$$
\|e^{-\eps t}\partial_t \th\|^2_{L^2(0,\iy;L^2(0,\pi))}
+\|e^{-\eps t} \th \|_{L^2((0,\iy);\cH_2)}\prec
\|e^{-\eps t}\partial_t f \|_{L^2((0,\iy),\cH_1)}
+\| \xi \|_{ \cH_2}.
$$
}

In the right hand side we see, roughly speaking, the $L^2$norm of
$\partial_t \partial_x f$.
We can slightly  strengthen this estimate with the $L^2$norm of
$\partial_t  f$ in the right hand side. It is a sharp result in the sense that it close
to the estimate \eqref{regwe2}

\BE{theorem}{\LB{strong}}{
Let \eqref{C1} is fulfilled,
$f_t\in \fH_{0,\gm}$, $f(x,0)=0$,
and $\xi\in \cH_s$. Then
$$
\|e^{-\eps t}\partial_t \th\|^2_{L^2(0,\iy;L^2(0,\pi))}
+\|e^{-\eps t} \th \|_{L^2((0,\iy);\cH_2)}\prec
\|e^{-\eps t}\partial_t f \|_{L^2((0,\iy),\cH_0)}
+\| \xi \|_{ \cH_2}.
$$
}
\proof

(i) Let $F=0$. Then
$$
\cL[\th_n']=z\cL [\th_n]-\th_n(0)=
\l(\frac{z}{z+n^2K(z)}-1  \r)\xi_n=-\frac{n^2K(z)}{z+n^2K(z)}\xi_n.
$$
$K(z)$ decreases in infinity,
$$
|K(1+iy)|\prec \frac1{1+|y|},
$$

Therefore we can use
$$
|\cL \th_n'|\prec n\frac{|K|}{|G_n|}|\xi_n|\prec n|\xi_n|,
$$
what gives by \eqref{2.2} the estimate
\EQ{f0}{
\| e^{-\eps t}\th_t \|\prec \| \xi \|_{H_{s+2}}.
}

(ii). Let $\xi=0$. Then, taking into account $f(x,0)=0$, we obtain
$$
\cL(\th_n')=\frac{zF_n(z)}{z+n^2K(z)}.
$$
By \eqref{2.1a} we have
$$
\|\cL(\th_n')\| \prec \|zF_n(z)\|.
$$
Then
\EQn{
\| e^{-\gm t}\th_n' \|\prec \|e^{-\gm t} f_n' \|_{L^2(0,\iy)}.
}
and
\EQ{xi0}{
\| e^{-\gm t}\th_t \|\prec \|e^{-\gm t} f' \|_{L^2(0,\iy;\cH_s)}.
}
The theorem is proved.

\textbf{Acknowledgements}

The author  grateful to Prof. V. V. Vlasov for the fruitful discussions.


\begin{thebibliography}{1}

\bibitem{Sham} A. A. Gavrikov, S.A. Ivanov, D.Yu. Knyazkov,
V.A. Samarain , A.S. Shamaev, V. V. Vlasov,
{\em Spectral properties of composite medaia},
Contemporary Problems of Mathematic and Mechanic, v.1, 2009, 142-159 (Russian).

\bibitem{GuPip} M. E. Gurtin, A. C. Pipkin
{\em A general theory of heat conduction with finite wave speeds}.
Archive for Rational Mechanics and Analysis 1968; 32:113-126.

\bibitem{CMD} C. M. Dafermos,  {\em Asymptotic stability in viscoelasticity},
  Arch. Rational Mech. Anal., 37 (1970), 297-308.

\bibitem{P05}
L. Pandolfi,{\em The controllability of the Gurtin-Pipkin equation:
a cosine operator approach}.
Appl. Math. Optim. 52 (2005), no. 2, 143--165.


\bibitem{EMV} F.M.  Vegni.
{\em Dissipativity of a consedved phase field systems with memory},
Discrete and contiouns dynamical systems ,
Volume 9, Number 4, July 2003.


\bibitem{VW}
V. V. Vlasov, J. Wu,
{\em Solvability and Spectral Analysis of Abstract Hyperbolic Equations with
Delay}.
Funct. Differ. Equ. 16 (2009), no. 4, 751-768.

\bibitem{I}Ivanov S.A.,
'Wave type' spectrum of the Gurtin-Pipkin equation of the second order,
arXiv; arxiv.org/abs/1002.2831, 8 p.


\bibitem{RSV} Vlasov V.V., Rautian N.A., Shamaev A.S., {\em
Spectral analysis and correct solvability of abstract integrodifferential
equations arising in thermophysics and acoustics}, Journal of Mathematical Sciences
April 2013, Volume 190, Issue 1, pp 34-65.

\bibitem{IE} Ivanov S.A., EremenkoA.,
{\em Spectra of the Gurtin-Pipkin type equations}, SIAM J. Math. Anal. 43, pp. 2296-2306.

\bibitem{Mikh} V. Mikhailov, {\em  Partial Differential Equations}Mir Publishers, Moscow, Russia, 1978.

\end{thebibliography}
\end{document}